\documentclass[12 pt]{article}
\usepackage{amsmath}
\usepackage{amssymb}
\usepackage{geometry}
\newtheorem {thm} {Theorem} [section]
\newtheorem {cly} {Corollary} [section]
\newtheorem{defn} {Definition} [section]
\newtheorem{eg} {Example} [section]
\newtheorem{rmrk} {Remark} [section]
\date{}
\begin{document}
\title{FG-coupled fixed point theorems for contractive and generalized quasi-contractive mappings}
\author{Deepa Karichery and Shaini Pulickakunnel}
\maketitle
\section*{Abstract}
In this paper we prove FG-coupled fixed point theorems for different contractive mappings and 
generalized quasi- contractive mappings in partially ordered complete metric spaces. We prove the existence of 
FG-coupled fixed points of continuous as well as discontinuous mappings. We give some examples to illustrate the results.\\

\hspace*{-.6cm}\textbf{Key words:} FG-coupled fixed point, mixed monotone property, partially ordered complete metric space, contractive mapping, quasi-contactive mapping.\\
\textbf{MSC(2010):} 47H10, 54F05
\section{Introduction and Preliminaries}
Coupled fixed point problems become a new trend in non-linear analysis as a generalization of fixed point theory. 
In 1987 Dajun Guo and V. Lakshmikantham \cite{6} introduced the concept of coupled fixed point and proved some 
coupled fixed point theorems for mixed monotone mappings in cone metric spaces. Later in 2006 
Gnana Bhaskar and Lakshmikantham \cite{5} introduced the concept of coupled fixed point and mixed monotone property 
for contractive mappings in partially ordered complete metric spaces. Thereafter by changing the spaces and using 
different contractions several authors have proved various coupled fixed point theorems \cite{1,2,7,9}. \\

Recently E. Prajisha and P. Shaini \cite{8} defined FG-coupled fixed point on product of two spaces as a generalization 
of coupled fixed point and some FG-fixed point theorems have been proved. In this paper we prove FG- coupled fixed point 
theorems for contractive type mappings on  partially ordered complete metric spaces. Our results generalizes several 
fixed point theorems in literature. Now we recall some definitions.
\begin{defn}[\cite{5}]\normalfont
An element $(x,y)\in X\times X$ is said to be a coupled fixed point of the map $F:X\times X\rightarrow X$ if $F(x,y)=x$ and $F(y,x)=y$. 
\end{defn}
\begin{defn}[\cite{8}]\normalfont
Let $(X,d_X,\leq_{P_1})$ and $(Y,d_Y,\leq_{P_2})$ be two partially ordered metric spaces and $F:X\times Y\rightarrow X$ and $G:Y\times X\rightarrow Y$. We say that F and G have mixed monotone property if for any $x,y \in X$\\
$x_1, x_2 \in X,\ \ x_1\leq_{P_1} x_2\Rightarrow F(x_1,y)\leq_{P_1} F(x_2,y)$ and $G(y,x_1)\geq_{P_2}G(y,x_2)$\\
$y_1,y_2 \in Y,\ \ y_1\leq_{P_2}y_2\Rightarrow F(x,y_1)\geq_{P_1}F(x,y_2)$ and $G(y_1,x)\leq_{P_2}G(y_2,x)$.
\end{defn}
\begin{defn}[\cite{8}]\normalfont
An element $(x,y)\in X\times Y$ is said to be FG- coupled fixed point if $F(x,y)=x$ and $G(y,x)=y$.
\end{defn}
\textbf{Notes:}\\
1) If $(x,y)\in X\times Y$ is an FG- coupled fixed point then $(y,x)\in Y\times X$ is GF- coupled\\ \hspace*{.3cm} fixed point.\\
2) The metric $d$ on $X\times Y$ is defined by $d((x,y),(u,v))=d_X(x,u)+d_Y(y,v)$ for all\\ \hspace*{.3cm} $(x,y),(u,v)\in X\times Y$.\\
3) Partial order $\leq$ on $X\times Y$ is defined by for any $(x,y), (u,v)\in X\times Y; (u,v)\leq (x,y)\\ \hspace*{.4cm}\Leftrightarrow x\geq_{P_1}u,\ y\leq_{P_2} v$. \\
4) $F^{n+1}(x,y)=F(F^n(x,y),G^n(y,x))$  and $G^{n+1}(y,x)=G(G^n(y,x),F^n(x,y))$ for every \hspace*{.5cm}$n\in \mathbb{N}$ and $(x,y)\in X\times Y$.\\

\section{FG-coupled fixed point theorems for contraction mappings}
\begin{thm}\label{thm1}
Let $(X,d_X,\leq_{P_1})$ and $(Y,d_Y,\leq_{P_2})$ be two partially ordered complete metric spaces and $F:X\times Y\rightarrow X$, $G:Y\times X\rightarrow Y$ be two continuous mappings having the mixed monotone property. Assume that there exist $k,l,m,n\in [0,1); k+l<1$ and $m+n<1$ with
\begin{equation}\label{eq1}
d_X(F(x,y),F(u,v))\leq k\ d_X(x,u)+l\ d_Y(y,v)\ \forall x\geq_{P_1}u,\ y\leq_{P_2}v
\end{equation} 
\begin{equation}\label{eq2}
d_Y(G(y,x),G(v,u))\leq m\ d_Y(y,v)+n\ d_X(x,u)\ \forall x\leq_{P_1}u,\ y\geq_{P_2}v.
\end{equation}
If there exist $x_0\leq_{P_1} F(x_0, y_0)$ and $y_0\geq_{P_2} G(y_0, x_0)$, then there exist $(x,y)\in X\times Y$ such that $x=F(x,y)$ and $y=G(y,x)$.
\end{thm}   
\textbf{Proof: } Given $x_0\leq_{P_1} F(x_0, y_0)=x_1$ (say) and $y_0\geq_{P_2} G(y_0, x_0)=y_1$ (say).\\
Define $x_{n+1}=F(x_n,y_n)$ and $y_{n+1}= G(y_n,x_n)$ for $n=1,2,3...$\\
Then we get $F^{n+1}(x_0,y_0)=x_{n+1}$ and $G^{n+1}(y_0,x_0)=y_{n+1}$. \\
Using mathematical induction and mixed monotone property of F and G we prove that $\lbrace x_n \rbrace$ is increasing in X and $\lbrace y_n \rbrace$ is decreasing in Y. For,\\
Given $x_0\leq_{P_1} x_1$ and $y_0\geq_{P_2} y_1$. Claim that $x_n \leq_{P_1} x_{n+1}\ \  and \ \ y_n\geq_{P_2} y_{n+1}\ \ \forall n\in \mathbb{N}.$\\
For $n=1$, $x_2=F(x_1,y_1)\geq_{P_1} F(x_0,y_1)\geq_{P_1} F(x_0, y_0)=x_1$ and\\
$y_2=G(y_1,x_1)\leq_{P_2} G(y_1,x_0)\leq_{P_2} G(y_0,x_0)= y_1$.\\
Assume that the result is true for n=m.
ie, 
$x_{m+1}\geq_{P_1} x_m\ \ and\ \ y_m\geq_{P_2} y_{m+1}$.\\
Now consider,\\ 
$x_{m+2}=F(x_{m+1},y_{m+1})\geq_{P_1} F(x_m,y_{m+1})\geq_{P_1} F(x_m,y_m)=x_{m+1}$\\
$y_{m+2}=G(y_{m+1},x_{m+1})\leq_{P_2} G(y_{m+1},x_m)\leq_{P_2} G(y_m,x_m)=y_{m+1}$\\
So the result is true for $\forall n\in \mathbb{N}$.\\
Now, for proving the sequences $\lbrace x_n \rbrace$ and $\lbrace y_n \rbrace$ are Cauchy, we consider two cases.\\
Case 1: $m+n\leq k+l$\\
Claim that, for $j\in \mathbb{N}$, \begin{equation}\label{eq3}
d_X(F^{j+1}(x_0,y_0),F^j(x_0,y_0))\leq (k+l)^j[d_X(x_1,x_0)+d_Y(y_1,y_0)]
\end{equation}
\begin{equation}\label{eq4}
d_Y(G^{j+1}(y_0,x_0),G^j(y_0,x_0))\leq (k+l)^j[d_Y(y_1,y_0)+d_X(x_1,x_0)]
\end{equation}
The proof is by mathematical induction.\\
For $j=1$, consider 
\begin{align*}
d_X(F^2(x_0,y_0),F(x_0,y_0))=&d_X(F(F(x_0,y_0),G(y_0,x_0)),F(x_0,y_0))\\
\leq&  k\ d_X(F(x_0,y_0),x_0)+ l\ d_Y(G(y_0,x_0),y_0)\\
\leq& (k+l)\ [d_X(x_1,x_0)+ d_Y(y_1,y_0)]
\end{align*}
Similarly we prove that $d_Y(G^2(y_0,x_0),G(y_0,x_0))\leq \ (k+l)[d_X(x_1,x_0)+d_Y(y_1,y_0)]$. \\
ie, the result is true for $j=1$.\\
Now assume that the claim is true for $j\leq t$, and prove for $j=t+1$. \\
Consider,\\ 
$d_X(F^{t+2}(x_0,y_0),F^{t+1}(x_0,y_0))\\ 
~~~~~~~~~~~= d_X(F(F^{t+1}(x_0,y_0),G^{t+1}(y_0,x_0)),F(F^{t}(x_0,y_0),G^{t}(y_0,x_0)))\\
~~~~~~~~~~~\leq  k\ d_X(F^{t+1}(x_0,y_0),F^t(x_0,y_0))+ l\ d_Y(G^{t+1}(y_0,x_0),G^t(y_0,x_0))\\
~~~~~~~~~~~\leq  k(k+l)^t[d_X(x_1,x_0)+d_Y(y_1,y_0)]+ l\ (k+l)^t [d_X(x_1,x_0)+d_Y(y_1,y_0)]\\
~~~~~~~~~~~ = (k+l)^{t+1}[d_X(x_1,x_0)+ d_Y(y_1,y_0)]$\\
Similarly we get $d_Y(G^{t+2}(y_0,x_0),G^{t+1}(y_0,x_0))\leq (k+l)^{t+1}[d_X(x_1,x_0)+d_Y(y_1,y_0)]$.\\
Thus the claim is true for all $j\in \mathbb{N}$.\\
Now we prove that $\lbrace F^j(x_0,y_0)\rbrace$ and $\lbrace G^j(y_0,x_0)\rbrace$ are Cauchy sequences in X and Y respectively.\\
For $t>j$, consider\\
$d_X(F^j(x_0,y_0),F^t(x_0,y_0))\\
~~~~~~~~~~~\leq  d_X(F^j(x_0,y_0),F^{j+1}(x_0,y_0))+d_X(F^{j+1}(x_0,y_0),F^{j+2}(x_0,y_0))+...\\
~~~~~~~~~~~~~~~+d_X(F^{t-1}(x_0,y_0),F^t(x_0,y_0))\\
~~~~~~~~~~~\leq \Big[(k+l)^j+(k+l)^{j+1}+...+(k+l)^{t-1}\Big][d_X(x_1,x_0)+d_Y(y_1,y_0)]\\
~~~~~~~~~~~\leq \dfrac{{\delta_1 }^j}{1-\delta_1}\Big[d_X(x_1,x_0)+d_Y(y_1,y_0)\Big]\ \text{where}\ \delta_1 = k+l <1\\
~~~~~~~~~~~\rightarrow 0\ \text{as}\ \ j\rightarrow \infty$\\
ie,  $\lbrace F^j(x_0,y_0)\rbrace$ is a Cauchy sequence in X.\\
Similarly we can prove that $\lbrace G^j(y_0,x_0)\rbrace$ is a Cauchy sequence in Y.\\
Case 2: $k+l < m+n$.\\
Now we claim that \begin{equation}\label{eq22}
d_X(F^{j+1}(x_0,y_0),F^j(x_0,y_0))< (m+n)^j[d_X(x_1,x_0)+d_Y(y_1,y_0)]
\end{equation}
\begin{equation}\label{eq23}
d_Y(G^{j+1}(y_0,x_0),G^j(y_0,x_0))< (m+n)^j[d_Y(y_1,y_0)+d_X(x_1,x_0)]
\end{equation}
For $j=1$
\begin{align*}
d_X(F^2(x_0,y_0),F(x_0,y_0))&=d_X(F(F(x_0,y_0),G(y_0,x_0)),F(x_0,y_0))\\
&\leq  k\ d_X(F(x_0,y_0),x_0)+ l\ d_Y(G(y_0,x_0),y_0)\\
&\leq  (k+l)\ [d_X(x_1,x_0)+ d_Y(y_1,y_0)]\\
&<  (m+n)\ [d_X(x_1,x_0)+ d_Y(y_1,y_0)]
\end{align*}
Similarly we can get $d_Y(G^2(y_0,x_0),G(y_0,x_0))< (m+n)[d_Y(y_1,y_0)+d_X(x_1,x_0)]$.\\
Assume that the claim is true for $j\leq t$. Now we prove the claim for $j=t+1$.\\
Consider,\\
$d_X(F^{t+2}(x_0,y_0),F^{t+1}(x_0,y_0))\\
~~~~~~~~~= d_X(F(F^{t+1}(x_0,y_0),G^{t+1}(y_0,x_0)),F(F^{t}(x_0,y_0),G^{t}(y_0,x_0)))\\
~~~~~~~~~\leq k\ d_X(F^{t+1}(x_0,y_0),F^t(x_0,y_0))+ l\ d_Y(G^{t+1}(y_0,x_0),G^t(y_0,x_0))\\
~~~~~~~~~<  k(m+n)^t[d_X(x_1,x_0)+d_Y(y_1,y_0)]+ l\ (m+n)^t [d_X(x_1,x_0)+d_Y(y_1,y_0)]\\
~~~~~~~~~<  (m+n)^{t+1}[d_X(x_1,x_0)+ d_Y(y_1,y_0)]$\\
Similarly we can prove that\\ $d_Y(G^{t+2}(y_0,x_0),G^{t+1}(y_0,x_0))< (m+n)^{t+1}[d_Y(y_1,y_0)+d_X(x_1,x_0)]$.
Hence the claim is true for all $j\in \mathbb{N}$. \\ 
Now we prove that $\lbrace F^j(x_0,y_0)\rbrace$ and $\lbrace G^j(y_0,x_0)\rbrace $ are Cauchy sequences in X and Y respectively.\\
For $t>j$, consider\\ 
$d_X(F^t(x_0,y_0),F^j(x_0,y_0))\\
~~~~~~~~~~~\leq  d_X(F^j(x_0,y_0),F^{j+1}(x_0,y_0))+ d_X(F^{j+1}(x_0,y_0),F^{j+2}(x_0,y_0))+...\\
~~~~~~~~~~~~~~~+d_X(F^{t-1}(x_0,y_0),F^{t}(x_0,y_0))\\
~~~~~~~~~~~<  \Big[(m+n)^j+(m+n)^{j+1}+...+ (m+n)^{t-1}\Big]\ [d_X(x_1,x_0)+d_Y(y_1,y_0)]\\
~~~~~~~~~~~\leq  \dfrac{\delta_2 ^j}{1-\delta_2}[d_X(x_1,x_0)+d_Y(y_1,y_0)]\ \text{where}\ \delta_2=m+n<1\\
~~~~~~~~~~~\rightarrow 0\ \text{as}\ j\rightarrow \infty$\\
So, $\lbrace F^j(x_0,y_0)\rbrace $ is a Cauchy sequence in X.\\
Similarly we can prove that $\lbrace G^j(y_0,x_0)\rbrace$ is a Cauchy sequence in Y.\\
Since X and Y are complete metric spaces, we have $\lim_{j\rightarrow \infty}F^j(x_0,y_0)=x$ and $\lim_{j\rightarrow \infty}G^j(y_0,x_0)=y$, where $x\in X$ and $y\in Y$.\\
Now we prove that $F(x,y)=x$ and $G(y,x)=y$.\\
Consider,
\begin{align*}
d_X(F(x,y),x)&= \lim_{j\rightarrow \infty}d_X(F(F^j(x_0,y_0),G^j(y_0,x_0)),F^j(x_0,y_0))\\
&=\lim_{j\rightarrow \infty} d_X(F^{j+1}(x_0,y_0),F^j(x_0,y_0))\\
&= 0 
\end{align*} 
Therefore $F(x,y)=x$. Similarly we can prove that $G(y,x)=y$.\hspace*{.25cm} $\square$
\begin{eg}\normalfont Let $X=(-\infty, 0]$ and $Y=[0,\infty)$ with usual order and usual metric. Define $F: X\times Y\rightarrow X$ by $F(x,y)=\dfrac{x}{3}-\dfrac{y}{4}$ and $G: Y\times X\rightarrow Y$ by $G(y,x)=\dfrac{y}{8}-\dfrac{x}{6}$. Then it is easy to check that F satisfies (\ref{eq1}) with $k=\dfrac{1}{3}, l=\dfrac{1}{4}$ and G satisfies (\ref{eq2}) with $m=\dfrac{1}{8}, n=\dfrac{1}{6}$ and $(0,0)$ is the FG-coupled fixed point.
\end{eg}
\begin{cly}{\upshape\lbrack 5, Theorem 2.1\rbrack}
Let $(X,\leq)$ be a partially ordered set and suppose there is a metric d on X such that $(X, d)$ is a complete metric
space. Let $F : X \times X \longrightarrow X $ be a continuous mapping having the mixed monotone property on X. Assume that there exists
$k\in[0,1)$ with
\begin{equation}\label{eq21}
d(F(x,y),F(u,v))\leq \dfrac{k}{2}\ [d(x,u)+d(y,v)]\ \ \forall x\geq u, y\leq v.
\end{equation}
If there exist $x_0,y_0\in X$ such that $x_0\leq F(x_0,y_0)$ and $y_0\geq F(y_0,x_0)$. Then there exist $x,y\in X$ such that $x=F(x,y)$ and $y=F(y,x)$.
\end{cly}
\textbf{Proof}: Taking $X=Y$, $F=G$ and $k=l=m=n=\dfrac{k}{2}$ in Theorem \ref{thm1} we obtain the result.\hspace*{.25cm} $\square$
\begin{rmrk}\normalfont It can be shown that FG- coupled fixed point is unique provided that for every $(x,y),(x^*,y^*)\in X\times Y$, there exist a $(z_1,z_2)\in X\times Y$ that is comparable to $(x,y)$ and $(x^*,y^*)$. The result is proved in the following theorem.
\end{rmrk}
\begin{thm}\label{thm9}
Let $(X,d_X,\leq_{P_1})$ and $(Y,d_Y,\leq_{P_2})$ be two partially ordered complete metric spaces and $F:X\times Y\rightarrow X$, $G:Y\times X\rightarrow Y$ be two continuous mappings having the mixed monotone property. For every $(x,y),(x^*,y^*)\in X\times Y$ there exist a $(z_1,z_2)\in X\times Y$ that is comparable to $(x,y)$ and $(x^*,y^*)$. Assume that there exist $k,l,m,n\in [0,1);\ k+l<1$ and $m+n<1$ with
\begin{equation*}
d_X(F(x,y),F(u,v))\leq k\ d_X(x,u)+l\ d_Y(y,v)\ \ \ \forall x\geq_{P_1}u,\ y\leq_{P_2}v
\end{equation*} 
\begin{equation*}
d_Y(G(y,x),G(v,u))\leq m\ d_Y(y,v)+n\ d_X(x,u)\ \forall x\leq_{P_1}u,\ y\geq_{P_2}v.
\end{equation*}
If there exist $x_0\leq_{P_1} F(x_0, y_0)$ and $y_0\geq_{P_2} G(y_0, x_0)$, then there exist a unique FG- coupled fixed point.
\end{thm}
\textbf{Proof:} The existence of FG- coupled fixed point is followed by the proof of Theorem \ref{thm1}. Now we prove that if $(x^*,y^*)$ is another FG- coupled fixed point then $d((x,y),(x^*,y^*))=0$ where $x=\lim_{j\rightarrow \infty}F^j(x_0,y_0)$ and $y=\lim_{j\rightarrow \infty}G^j(y_0,x_0)$. \\
Case 1: If $(x,y)$ is comparable to $(x^*,y^*)$ with respect to the ordering in $X\times Y$, then $(F^j(x,y),G^j(y,x))=(x,y)$ is comparable to $(F^j(x^*,y^*),G^j(y^*,x^*))=(x^*,y^*)$ for every $j=1,2,3...$\\
If $m+n\leq k+l$, consider, \\
$d((x,y),(x^*,y^*))=d_X(x,x^*)+d_Y(y,y^*)\\
\hspace*{3cm}=d_X(F^j(x,y),F^j(x^*,y^*))+d_Y(G^j(y,x),G^j(y^*,x^*))\\
\hspace*{3cm}\leq 2^j(k+l)^j[d_X(x,x^*)+d_Y(y,y^*)]\\
\hspace*{3cm}\rightarrow 0\ \text{as}\ j\rightarrow\infty $.\\
Which implies that $d((x,y),(x^*,y^*))=0$.\\
Similarly for $k+l< m+n$ we get\\
$d((x,y),(x^*,y^*))< 2^j (m+n)^j[d_X(x,x^*)+d_Y(y,y^*)]\\
\hspace*{3cm}\rightarrow 0\ \text{as}\ j\rightarrow\infty $.\\
Therefore $d((x,y),(x^*,y^*))=0$.\\
Case 2: If $(x,y)$ is not comparable to $(x^*,y^*)$, then there exist $(z_1,z_2)\in X\times Y$ such that which is comparable to $(x,y)$ and $(x^*,y^*)$.\\
Without loss of generality, consider $m+n\leq k+l$, then\\
$d((x,y),(x^*,y^*))\\
\hspace*{.5cm}=d((F^j(x,y),G^j(y,x)),(F^j(x^*,y^*),G^j(y^*,x^*)))\\
\hspace*{.5cm}\leq d((F^j(x,y),G^j(y,x)),(F^j(z_1,z_2),G^j(z_2,z_1)))+\\ \hspace*{1cm} d((F^j(z_1,z_2),G^j(z_2,z_1)),(F^j(x^*,y^*),G^j(y^*,x^*)))\\
\hspace*{.5cm}=d_X(F^j(x,y),F^j(z_1,z_2))+d_Y(G^j(y,x),G^j(z_2,z_1))+ d_X(F^j(z_1,z_2),F^j(x^*,y^*))+\\ \hspace*{1cm} d_Y(G^j(z_2,z_1),G^j(y^*,x^*))\\
\hspace*{.5cm}\leq 2^j(k+l)^j\lbrace[d_X(x,z_1)+d_Y(y,z_2)]+[d_X(z_1,x^*)+d_Y(z_2,y^*)]\rbrace\\
\hspace*{.5cm}\rightarrow 0\ \text{as}\ j\rightarrow\infty $.\\
 Therefore $d((x,y),(x^*,y^*))=0$.\hspace*{.25cm} $\square$\\

We can replace the continuity of F and G by other conditions to get the existence of FG- coupled fixed point, as shown in the following theorem.
\begin{thm}\label{thm7}
Let $(X,d_X,\leq_{P_1})$ and $(Y,d_Y,\leq_{P_2})$ be two partially ordered complete metric spaces and $F:X\times Y\rightarrow X$, $G:Y\times X\rightarrow Y$ be two mappings satisfying the mixed monotone property. Assume that X and Y having the following property
\begin{enumerate}
\item[(i)] If a non-decreasing sequence $\lbrace x_n \rbrace\rightarrow x$ then $x_n\leq_{P_1} x\ \forall n$
\item[(ii)] If a non-increasing sequence $\lbrace y_n \rbrace\rightarrow y$ then $y\leq_{P_2} y_n\ \forall n$. 
\end{enumerate}
Also assume that there exist $k,l,m,n\in [0,1)$ such that $k+l<1$, $m+n<1$ with
\begin{equation}\label{eq34}
d_X(F(x,y),F(u,v))\leq k\ d_X(x,u)+l\ d_Y(y,v)\ \ \ \forall x\geq_{P_1}u,\ y\leq_{P_2}v
\end{equation} 
\begin{equation}\label{eq35}
d_Y(G(y,x),G(v,u))\leq m\ d_Y(y,v)+n\ d_X(x,u)\ \ \ \forall x\leq_{P_1}u,\ y\geq_{P_2}v.
\end{equation}
If there exist $x_0\leq_{P_1} F(x_0, y_0)$ and $y_0\geq_{P_2} G(y_0, x_0)$, then there exist $(x,y)\in X\times Y$ such that $x=F(x,y)$ and $y=G(y,x)$.
\end{thm}
\textbf{Proof:} Following as in the proof of Theorem \ref{thm1}, we get $\lim_{j\rightarrow \infty}F^j(x_0,y_0)=x$ and $\lim_{j\rightarrow \infty}G^j(y_0,x_0)=y$. \\
Now we have,
\begin{align*}
d_X(F(x,y),x)&\leq d_X(F(x,y),F^{j+1}(x_0,y_0))+ d_X(F^{j+1}(x_0,y_0),x)\\
&= d_X(F(x,y),F(F^j(x_0,y_0),G^j(y_0,x_0)))+ d_X(F^{j+1}(x_0,y_0),x)\\
&\leq  k\ d_X(x,F^j(x_0,y_0))+ l\ d_Y(y,G^j(y_0,x_0))+  d_X(F^{j+1}(x_0,y_0),x)\\
&\rightarrow 0\ \text{as}\ j\rightarrow \infty
\end{align*}
Therefore  $x=F(x,y)$.\\
Similarly using (\ref{eq35}) and $\lim_{j\rightarrow \infty}G^j(y_0,x_0)=y$, we get $y=G(y,x)$.\hspace*{.25cm} $\square$

\begin{cly} {\upshape\lbrack 5, Theorem 2.2\rbrack}
Let $(X,\leq)$ be a partially ordered set and suppose there is a metric d on X such that $(X, d)$ is a complete metric
space. Assume that X has the following property
\begin{enumerate}
\item[(i)] If a non-decreasing sequence $\lbrace x_n \rbrace\rightarrow x$ then $x_n\leq_{P_1} x\ \forall n$
\item[(ii)] If a non-increasing sequence $\lbrace y_n \rbrace\rightarrow y$ then $y\leq_{P_2} y_n\ \forall n$. 
\end{enumerate}
Let $F : X \times X \longrightarrow X $ be a mapping having the mixed monotone property on X. Assume that there exist  $k\in [0,1)$ such that \begin{equation*}
d(F(x,y),F(u,v))\leq \dfrac{k}{2}\ [d(x,u)+d(y,v)]\ \forall x\geq u, y\leq v.
\end{equation*} If there exist $x_0,y_0\in X$ such that $x_0\leq F(x_0,y_0)$ and $y_0\geq F(y_0,x_0)$, then there exist $x,y\in X$ such that $x=F(x,y)$ and $y=F(y,x)$.
\end{cly}
\textbf{Proof:} Taking $X=Y$, $F=G$ and $k=l=m=n=\dfrac{k}{2}$ in Theorem \ref{thm7}, we get the result.\hspace*{.25cm} $\square$
\begin{rmrk}\normalfont
In Theorem \ref{thm7} if add the condition: for every $(x,y),(x^*,y^*)\in X\times Y$, there exist a $(z_1,z_2)\in X\times Y$ that is comparable to both $(x,y)$ and $(x^*,y^*)$, we get unique FG- coupled fixed point.
\end{rmrk}
\begin{rmrk}\normalfont
If we take $k=l=\dfrac{a}{2}$ and $m=n=\dfrac{b}{2}$ where $a,b\in [0,1)$ with $a+b<1$ in Theorems \ref{thm1}, \ref{thm9} and \ref{thm7}, we get Theorems 2.1, 2.2 and 2.3 respectively of \cite{8}.
\end{rmrk}
\begin{rmrk}\normalfont
If we take $k=m$ and $l=n$ in Theorems \ref{thm1}, \ref{thm9} and \ref{thm7}, we get Theorems 2.4, 2.5 and 2.6 respectively of \cite{8}.
\end{rmrk}
\begin{thm}\label{thm2}
Let $(X,d_X,\leq_{P_1}), (Y,d_Y,\leq_{P_2})$ be two partially ordered complete metric spaces. Let $F:X\times Y\rightarrow X$ and $G:Y\times X\rightarrow Y$ be two continuous functions having the mixed monotone property. Assume that there exist $k,l,m,n\in \Big[0,\dfrac{1}{2}\Big)$ satisfying\begin{equation}\label{eq5}
d_X(F(x,y),F(u,v))\leq k\ d_X(x,F(x,y))+l\ d_X(u, F(u,v)); \forall x\geq_{P_1}u,\ y\leq_{P_2}v
\end{equation}
\begin{equation}\label{eq6}
d_Y(G(y,x),G(v,u))\leq m\ d_Y(y,G(y,x))+n\ d_Y(v,G(v,u)); \forall x\leq_{P_1}u, y\geq_{P_2}v.
\end{equation}
If there exist $x_0\in X, y_0\in Y$ satisfying  $x_0\leq_{P_1} F(x_0, y_0)$ and $y_0\geq_{P_2} G(y_0, x_0)$ then there exist $x\in X, y\in Y$ such that $x=F(x,y)$ and $y=G(y,x)$.
\end{thm}
\textbf{ Proof:} Following as in Theorem \ref{thm1} we can show that $\lbrace x_n \rbrace$ is increasing in X and $\lbrace y_n \rbrace$ is decreasing in Y.\\
Using inequalities (\ref{eq5}) and (\ref{eq6}) we get
\begin{align*}
d_X(x_{n+1},x_n)&= d_X(F(x_n,y_n),F(x_{n-1},y_{n-1}))\\
&\leq  k\ d_X(x_n, F(x_n,y_n))+ l\ d_X(x_{n-1},F(x_{n-1},y_{n-1}))\\
&= k\ d_X(x_n,x_{n+1})+ l\ d_X(x_{n-1},x_n)\\
\text{ie,}\ (1-k)\ d_X(x_{n+1},x_n)&\leq l\ d_X(x_{n-1},x_n)\\
\text{ie,}\ d_X(x_n,x_{n+1})&\leq  \dfrac{l}{1-k} d_X(x_{n-1},x_n)\\
&=\delta_1\ d_X(x_{n-1},x_n)\ \text{where}\ \ \delta_1= \dfrac{l}{1-k} < 1\\
&\leq \delta_1^2\ d_X(x_{n-2},x_{n-1})&&
\end{align*}
\begin{align*}
&\hspace*{1cm}\vdots &&\\
&\leq \delta_1^n\ d_X(x_0,x_1)
\end{align*}
Similarly we get $d_Y(y_{n+1},y_n)\leq {\delta_2}^n\  d_Y(y_1,y_0) $ where $\delta_2=\dfrac{m}{1-n}$\\ 
Consider $m>n$
\begin{align*}
d_X(x_m,x_n)&\leq d_X(x_m,x_{m-1})+d_X(x_{m-1},x_{m-2})+...+d_X(x_{n+1},x_n)\\
&\leq  {\delta_1}^{m-1}\ d_X(x_1,x_0)+ {\delta_1}^{m-2}\ d_X(x_1,x_0)+...+ {\delta_1}^n\ d_X(x_1,x_0)\\
&= {\delta_1}^n \Big(1+\delta_1 + ...+ {\delta_1}^{m-n-1}\Big)\ d_X(x_1,x_0)\\
&= \frac{{\delta_1}^n}{1-\delta_1}\ d_X(x_1,x_0)\\
&\rightarrow  0\ \text{as}\ n\rightarrow \infty.  
\end{align*}
Therefore $\lbrace F^n(x_0,y_0)\rbrace$ is a Cauchy sequence in X.\\
Similarly we can prove that $\lbrace G^n(y_0,x_0)\rbrace$ is a Cauchy sequence in Y.\\
Since by the completeness of X and Y, there exist $x\in X$ and $y\in Y$ such that $\lim _{n\rightarrow \infty}F^n(x_0,y_0)=x$ and $\lim _{n\rightarrow \infty}G^n(y_0,x_0)=y$.\\
As in the proof of Theorem \ref{thm1} by continuity of F and G we can prove that $F(x,y)=x$ and $G(y,x)=y$.\hspace*{.25cm} $\square$\\

By replacing the continuity of F and G by other conditions we obtain the following existence theorems of FG-coupled fixed point.
\begin{thm}\label{thm10}
Let $(X,d_X,\leq_{P_1})$ and $(Y,d_Y,\leq_{P_2})$ be two partially ordered complete metric spaces and $F:X\times Y\rightarrow X$, $G:Y\times X\rightarrow Y$ be two mappings having the mixed monotone property. Assume that X and Y satisfy the following property\begin{enumerate}
\item[(i)]. If a non-decreasing sequence $\lbrace x_n \rbrace\rightarrow x$ then $x_n\leq_{P_1} x\ \forall n$
\item[(ii)]. If a non-increasing sequence $\lbrace y_n \rbrace\rightarrow y$ then $y\leq_{P_2} y_n\ \forall n$ 
\end{enumerate}
Also assume that there exist $k,l,m,n\in \Big[0,\dfrac{1}{2}\Big)$ satisfying\begin{equation}\label{eq24}
d_X(F(x,y),F(u,v))\leq k\ d_X(x,F(x,y))+l\ d_X(u, F(u,v)); \forall x\geq_{P_1}u,\ y\leq_{P_2}v
\end{equation}
\begin{equation}\label{eq25}
d_Y(G(y,x),G(v,u))\leq m\ d_Y(y,G(y,x))+n\ d_Y(v,G(v,u)); \forall x\leq_{P_1}u,\ y\geq_{P_2}v.
\end{equation}
If there exist $x_0\in X,\ y_0\in Y$ satisfying  $x_0\leq_{P_1} F(x_0, y_0)$ and $y_0\geq_{P_2} G(y_0, x_0)$ then there exist $x\in X,\ y\in Y$ such that $x=F(x,y)$ and $y=G(y,x)$.
\end{thm}
\textbf{Proof}: Following as in the proof of Theorem \ref{thm2} we get $\lim_{n\rightarrow \infty}F^n(x_0,y_0)=x$ and $\lim_{n\rightarrow \infty}G^n(y_0,x_0)=y$. \\
Now we have
\begin{align*}
d_X(F(x,y),x)&\leq d_X(F(x,y),F^{n+1}(x_0,y_0))+ d_X(F^{n+1}(x_0,y_0),x)\\
&= d_X(F(x,y),F(F^n(x_0,y_0),G^n(y_0,x_0))+ d_X(F^{n+1}(x_0,y_0),x)\\
&\leq  k\ d_X(x,F(x,y))+ l\ d_X(F^n(x_0,y_0),F(F^n(x_0,y_0),G^n(y_0,x_0)))\\&\ \ \ \ + d_X(F^{n+1}(x_0,y_0),x)\ \ (\text{using (\ref{eq24})})
\end{align*}
ie, $d_X(F(x,y),x)\leq k\ d_X(x,F(x,y))$ as $n\rightarrow \infty$\\
This holds only when $d_X(F(x,y),x)=0$. Therefore we get $F(x,y)=x$.\\
Similarly using (\ref{eq25}) and $\lim_{n\rightarrow \infty}G^n(y_0,x_0)=y$ we can prove $y=G(y,x)$.\hspace*{.25cm} $\square$
\begin{rmrk}\normalfont
If we put $k=m$ and $l=n$ in Theorems \ref{thm2} and \ref{thm10}, we get Theorems 2.7 and 2.8 respectively of \cite{8}.
\end{rmrk}
\begin{thm}\label{thm3}
Let $(X,d_X,\leq_{P_1}), (Y,d_Y,\leq_{P_2})$ be two partially ordered complete metric spaces. Let $F:X\times Y\rightarrow X$ and $G:Y\times X\rightarrow Y$ be two continuous functions having the mixed monotone property. Assume that there exist $k,l,m,n\in \Big[0,\dfrac{1}{2}\Big)$ satisfying \begin{equation}\label{eq7}
d_X(F(x,y),F(u,v))\leq k\ d_X(x,F(u,v))+l\ d_X(u, F(x,y));\ \forall\ x\geq_{P_1}u,\ y\leq_{P_2}v
\end{equation}
\begin{equation}\label{eq8}
d_Y(G(y,x),G(v,u))\leq m\ d_Y(y,G(v,u))+n\ d_Y(v,G(y,x));\ \forall\ x\leq_{P_1}u,\ y\geq_{P_2}v.
\end{equation}
If there exist $x_0\in X, y_0\in Y$ satisfying  $x_0\leq_{P_1} F(x_0, y_0)$ and $y_0\geq_{P_2} G(y_0, x_0)$ then there exist $x\in X, y\in Y$ such that $x=F(x,y)$ and $y=G(y,x)$.
\end{thm}
\textbf{Proof:} As in Theorem \ref{thm1} we have $\lbrace x_n \rbrace$ is increasing in X and $\lbrace y_n \rbrace$ is decreasing in Y.\\
We have
\begin{align*}
d_X(x_{n+1},x_n)&= d_X(F(x_n,y_n),F(x_{n-1},y_{n-1}))\\
&\leq  k\ d_X(x_n, F(x_{n-1},y_{n-1}))+ l\ d_X(x_{n-1},F(x_n,y_n))\ \ (\text{Using (\ref{eq7})})\\
&= k\ d_X(x_n,x_n)+ l\ d_X(x_{n-1},x_{n+1})\\
&\leq l\ [d_X(x_{n-1},x_n)+d_X(x_n,x_{n+1})]\\
\text{ie,}\ \ d_X(x_n,x_{n+1})&\leq  \dfrac{l}{1-l}\ d_X(x_{n-1},x_n)\\
&=\delta_1\ d_X(x_{n-1},x_n)\ \ where\ \ \delta_1= \dfrac{l}{1-l} < 1\\
&\leq \delta_1^2\ d_X(x_{n-2},x_{n-1})\\
&\hspace*{1cm}\vdots &&\\
&\leq \delta_1^n\ d_X(x_0,x_1)
\end{align*}  
Similarly we get $d_Y(y_{n+1},y_n)\leq  {\delta_2}^n d_Y(y_1,y_0) $ where $\delta_2=\dfrac{m}{1-m}$\\
Now we prove that $\lbrace F^n(x_0,y_0) \rbrace$ and $\lbrace G^n(y_0,x_0) \rbrace$ are Cauchy sequences in X and Y respectively.\\
For $m>n$, \begin{align*}
d_X(x_m,x_n)&\leq d_X(x_m,x_{m-1})+d_X(x_{m-1},x_{m-2})+...+d_X(x_{n+1},x_n)\\
&\leq  {\delta_1}^{m-1}\ d_X(x_1,x_0)+ {\delta_1}^{m-2}\ d_X(x_1,x_0)+...+ {\delta_1}^n\ d_X(x_1,x_0)\\
&\leq \frac{{\delta_1}^n}{1-\delta_1}\ d_X(x_1,x_0)\\
&\rightarrow  0\ \text{as}\ n\rightarrow \infty.
\end{align*}
Therefore $\lbrace F^n(x_0,y_0) \rbrace$ is a Cauchy sequence in X.\\  
Similarly we can prove that $\lbrace G^n(y_0,x_0) \rbrace$ is a Cauchy sequence in Y.\\ 
By the completeness of X and Y, there exist $x\in X$ and $y\in Y$ such that $\lim _{n\rightarrow \infty}F^n(x_0,y_0)=x$ and $\lim _{n\rightarrow \infty} G^n(y_0,x_0)=y$.\\
As in the proof of Theorem \ref{thm1} we can show that $x=F(x,y)$ and $y=G(y,x)$. \hspace*{.25cm} $\square$\\
\begin{thm}\label{thm11}
Let $(X,d_X,\leq_{P_1})$ and $(Y,d_Y,\leq_{P_2})$ be two partially ordered complete metric spaces and $F:X\times Y\rightarrow X$, $G:Y\times X\rightarrow Y$ be two mappings having the mixed monotone property. Assume that X and Y satisfy the following property\begin{enumerate}
\item[(i)] If a non-decreasing sequence $\lbrace x_n \rbrace\rightarrow x$ then $x_n\leq_{P_1} x\ \forall n$
\item[(ii)] If a non-increasing sequence $\lbrace y_n \rbrace\rightarrow y$ then $y\leq_{P_2} y_n\ \forall n$. 
\end{enumerate}
Also assume that there exist $k,l,m,n\in \Big[0,\dfrac{1}{2}\Big)$ satisfying \begin{equation}\label{eq26}
d_X(F(x,y),F(u,v))\leq k\ d_X(x,F(u,v))+l\ d_X(u, F(x,y)); \forall x\geq_{P_1}u,\ y\leq_{P_2}v
\end{equation}
\begin{equation}\label{eq27}
d_Y(G(y,x),G(v,u))\leq m\ d_Y(y,G(v,u))+n\ d_Y(v,G(y,x)); \forall x\leq_{P_1}u,\ y\geq_{P_2}v.
\end{equation}
If there exist $x_0\in X,\ y_0\in Y$ satisfying  $x_0\leq_{P_1} F(x_0, y_0)$ and $y_0\geq_{P_2} G(y_0, x_0)$ then there exist $x\in X,\ y\in Y$ such that $x=F(x,y)$ and $y=G(y,x)$.
\end{thm}
\textbf{Proof:} Following as in the proof of Theorem \ref{thm3} we get $\lim_{n\rightarrow \infty}F^n(x_0,y_0)=x$ and $\lim_{n\rightarrow \infty}G^n(y_0,x_0)=y$.\\
Consider
\begin{align*}
d_X(F(x,y),x)&\leq d_X(F(x,y),F^{n+1}(x_0,y_0))+ d_X(F^{n+1}(x_0,y_0),x)\\
&= d_X(F(x,y),F(F^n(x_0,y_0),G^n(y_0,x_0)))+ d_X(F^{n+1}(x_0,y_0),x)\\
&\leq  k\ d_X(x,F((F^n(x_0,y_0),G^n(y_0,x_0)))+ l\ d_X(F^n(x_0,y_0),F(x,y))\\&\ \ \ + d_X(F^{n+1}(x_0,y_0),x)\\
&= k\ d_X(x,F^{n+1}(x_0,y_0))+ l\ d_X(F^n(x_0,y_0),F(x,y))\\&\ \ \ + d_X(F^{n+1}(x_0,y_0),x)
\end{align*}
ie, $d_X(F(x,y),x)\leq l\ d_X(x,F(x,y))$ as $n\rightarrow \infty$, which implies that $d_X(F(x,y),x)=0$. Therefore we get $F(x,y)=x$.\\
Similarly using (\ref{eq27}) and $\lim_{n\rightarrow \infty}G^n(y_0,x_0)=y$, we get $y=G(y,x)$. \hspace*{.25cm} $\square$
\begin{rmrk}\normalfont
If we put $k=m$ and $l=n$ in Theorems \ref{thm3} and \ref{thm11}, we get Theorems 2.9 and 2.10 respectively of \cite{8}.
\end{rmrk}
\begin{thm}\label{thm4}
Let $(X,d_X,\leq_{P_1}), (Y,d_Y,\leq_{P_2})$ be two partially ordered complete metric spaces. Let $F:X\times Y\rightarrow X$ and $G:Y\times X\rightarrow Y$ be two continuous functions having the mixed monotone property. Assume that there exist $a,b,c$ with $a+b+c<1$ satisfying\begin{equation}\label{eq9}
d_X(F(x,y),F(u,v))\leq a\ d_X(x,F(x,y))+b\ d_X(u, F(u,v))+c\ d_X(x,u); \forall x\geq_{P_1}u,\ y\leq_{P_2}v
\end{equation}
\begin{equation}\label{eq10}
d_Y(G(y,x),G(v,u))\leq a\ d_Y(y,G(y,x))+b\ d_Y(v,G(v,u))+c\ d_Y(y,v); \forall x\leq_{P_1}u,\ y\geq_{P_2}v.
\end{equation}
If there exist $x_0\in X, y_0\in Y$ satisfying  $x_0\leq_{P_1} F(x_0, y_0)$ and $y_0\geq_{P_2} G(y_0, x_0)$ then there exist $x\in X, y\in Y$ such that $x=F(x,y)$ and $y=G(y,x)$.
\end{thm}
\textbf{Proof:} Following as in Theorem \ref{thm1} we have $\lbrace x_n \rbrace$ is increasing in X and $\lbrace y_n \rbrace$ is decreasing in Y.\\ 
Now we claim that \begin{equation}\label{eq11}
d_X(F^{n+1}(x_0,y_0),F^n(x_0,y_0))\leq \Big(\dfrac{b+c}{1-a}\Big)^n\ d_X(x_0,x_1)
\end{equation}
\begin{equation}\label{eq12}
d_Y(G^{n+1}(y_0,x_0),G^n(y_0,x_0))\leq \Big(\frac{a+c}{1-b}\Big)^n\ d_Y(y_0,y_1)
\end{equation}
The proof is by mathematical induction with the help of (\ref{eq9}) and (\ref{eq10}).\\ For n=1, consider
\begin{align*}
d_X(F^2(x_0,y_0),F(x_0,y_0))&= d_X(F(F(x_0,y_0),G(y_0,x_0)),F(x_0,y_0))\\
&\leq  a\ d_X(F(x_0,y_0),F^2(x_0,y_0))+b\ d_X(x_0,F(x_0,y_0))+ \\&\ \ \ \ \ c\ d_X(F(x_0,y_0),x_0)\\
ie,\ d_X(F^2(x_0,y_0),F(x_0,y_0))&\leq  \dfrac{b+c}{1-a}\ d_X(x_0,x_1) 
\end{align*}
Thus the inequality (\ref{eq11}) is true for $n=1$.\\
Now assume that (\ref{eq11}) is true for $n\leq m$, and check for $n=m+1$.\\
Consider,\\
$d_X(F^{m+2}(x_0,y_0),F^{m+1}(x_0,y_0))\\
~~~~~~~~~~= d_X(F(F^{m+1}(x_0,y_0),G^{m+1}(y_0,x_0)),F(F^m(x_0,y_0),G^m(y_0,x_0)))\\
~~~~~~~~~~\leq  a\ d_X(F^{m+1}(x_0,y_0),F^{m+2}(x_0,y_0))+ b\ d_X(F^m(x_0,y_0),F^{m+1}(x_0,y_0))\\
~~~~~~~~~~~~~~~~~~~+ c\ d_X(F^{m+1}(x_0,y_0),F^m(x_0,y_0))\\
ie,\ d_X(F^{m+2}(x_0,y_0),F^{m+1}(x_0,y_0))\leq  \dfrac{b+c}{1-a}\ d_X(F^m(x_0,y_0),F^{m+1}(x_0,y_0))\\
\hspace*{6.2cm}\leq  \Big(\dfrac{b+c}{1-a}\Big)^{m+1}\ d_X(x_0,x_1)$\\
ie, the inequality (\ref{eq11}) is true for all $n\in \mathbb{N}$.\\
Similarly we can prove the inequality (\ref{eq12}).\\ 
For $m>n$, consider\\
$d_X(F^n(x_0,y_0),F^m(x_0,y_0))\\
~~~~~~~~\leq  d_X(F^n(x_0,y_0),F^{n+1}(x_0,y_0))+ d_X(F^{n+1}(x_0,y_0),F^{n+2}(x_0,y_0))+...\\
~~~~~~~~~~~+d_X(F^{m-1}(x_0,y_0),F^m(x_0,y_0))\\
~~~~~~~~\leq  \Bigg[\Big(\dfrac{b+c}{1-a}\Big)^n+\Big(\dfrac{b+c}{1-a}\Big)^{n+1}+...+\Big(\dfrac{b+c}{1-a}\Big)^{m-1}\Bigg]\ d_X(x_0,x_1)\\
~~~~~~~~\leq  \dfrac{{\delta_1}^n}{1-\delta_1}\ d_X(x_0,x_1)\ \text{where}\  \delta_1= \dfrac{b+c}{1-a}<1\\
~~~~~~~~\rightarrow 0\ \text{as}\ n\rightarrow \infty $\\
ie, $\lbrace F^n(x_0,y_0)\rbrace$ is a Cauchy sequence in X.\\
Similarly by using inequality (\ref{eq12}) we can prove that $\lbrace G^n(y_0,x_0)\rbrace$ is a Cauchy sequence in Y. \\
By the completeness of X and Y, there exist $x\in X$ and $y\in Y$ such that $\lim _{n\rightarrow \infty}F^n(x_0,y_0)=x$ and $\lim _{n\rightarrow \infty}G^n(y_0,x_0)=y$.\\
As in the proof of Theorem \ref{thm1}, using continuity of $F$ and $G$ we can prove that $F(x,y)=x$ and $G(y,x)=y$. \hspace*{.25cm} $\square$
\begin{thm}\label{thm12}
Let $(X,d_X,\leq_{P_1})$ and $(Y,d_Y, \leq_{P_2})$ be two partially ordered complete metric spaces and $F:X\times Y\rightarrow X$, $G:Y\times X\rightarrow Y$ be two mappings having the mixed monotone property. Assume that X and Y satisfy the following property\begin{enumerate}
\item[(i)] If a non-decreasing sequence $\lbrace x_n \rbrace\rightarrow x$ then $x_n\leq_{P_1} x\ \forall n$
\item[(ii)] If a non-increasing sequence $\lbrace y_n \rbrace\rightarrow y$ then $y\leq_{P_2} y_n\ \forall n$.
\end{enumerate}
Also asuume that there exist $a,b,c$ with $a+b+c<1$ satisfying\begin{equation}\label{eq28}
d_X(F(x,y),F(u,v))\leq a\ d_X(x,F(x,y))+b\ d_X(u, F(u,v))+c\ d_X(x,u); \forall x\geq_{p_1}u,\ y\leq_{p_2}v
\end{equation}
\begin{equation}\label{eq29}
d_Y(G(y,x),G(v,u))\leq a\ d_Y(y,G(y,x))+b\ d_Y(v,G(v,u))+c\ d_Y(y,v); \forall x\leq_{P_1}u,\ y\geq_{P_2}v.
\end{equation}
If there exist $x_0\in X, y_0\in Y$ satisfying  $x_0\leq_{P_1} F(x_0, y_0)$ and $y_0\geq_{P_2} G(y_0, x_0)$ then there exist $x\in X, y\in Y$ such that $x=F(x,y)$ and $y=G(y,x)$.
\end{thm}
\textbf{Proof}: Following as in the proof of Theorem \ref{thm4} we obtain $\lim_{n\rightarrow \infty}F^n(x_0,y_0)=x$ and $\lim_{n\rightarrow \infty}G^n(y_0,x_0)=y$.\\
We have
\begin{align*}
d_X(F(x,y),x)&\leq d_X(F(x,y),F^{n+1}(x_0,y_0))+ d_X(F^{n+1}(x_0,y_0),x)\\
&= d_X(F(x,y),F(F^n(x_0,y_0),G^n(y_0,x_0)))+ d_X(F^{n+1}(x_0,y_0),x)\\
&\leq  a\ d_X(x,F(x,y))+ b\ d_X(F^n(x_0,y_0),F(F^n(x_0,y_0),G^n(y_0,x_0)))+\\&\ \ \ \ c\ d_X(x,F^n(x_0,y_0))+ d_X(F^{n+1}(x_0,y_0),x)\\
&= a\ d_X(x,F(x,y))+ b\  d_X(F^n(x_0,y_0),F^{n+1}(x_0,y_0))\\&\ \ \ \ +c\ d_X(x,F^n(x_0,y_0))+ d_X(F^{n+1}(x_0,y_0),x)
\end{align*}
ie, $d_X(F(x,y),x)\leq a\ d_X(x,F(x,y))$ as $n\rightarrow \infty$, which implies that $d_X(F(x,y),x)=0$. Therefore $F(x,y)=x$.\\
Similarly using (\ref{eq29}) and $\lim_{n\rightarrow \infty}G^n(y_0,x_0)=y$ we get $y=G(y,x)$. \hspace*{.25cm} $\square$
\begin{rmrk}\normalfont
If we take $c=0$ in Theorems \ref{thm4} and \ref{thm12}, we get Theorems 2.7 and 2.8 respectively of \cite{8}.
\end{rmrk}
\begin{thm}\label{thm5}
Let $(X,d_X,\leq_{P_1}), (Y,d_Y,\leq_{P_2})$ be two partially ordered complete metric spaces. Let $F:X\times Y\rightarrow X$ and $G:Y\times X\rightarrow Y$ be two continuous functions having the mixed monotone property. Assume that there exist non-negative a,b,c satisfying\begin{equation}\label{eq13}
\begin{aligned}
d_X(F(x,y),F(u,v))\leq a\ d_X(x,F(u,v))+ &\ b\ d_X(u, F(x,y))+c\ d_X(x,u);\\& \forall x\geq_{P_1}u,\ y\leq_{P_2}v;\ 2b+c<1
\end{aligned}
\end{equation}
\begin{equation}\label{eq14}
\begin{aligned}
d_Y(G(y,x),G(v,u))\leq a\ d_Y(y,G(v,u))+ &\ b\ d_Y(v,G(y,x))+c\ d_Y(y,v);\\& \forall x\leq_{P_1}u,\ y\geq_{P2}v;\ 2a+c<1
\end{aligned}
\end{equation}
If there exist $x_0\in X, y_0\in Y$ satisfying  $x_0\leq_{P_1} F(x_0, y_0)$ and $y_0\geq_{P_2} G(y_0, x_0)$ then there exist $x\in X, y\in Y$ such that $x=F(x,y)$ and $y=G(y,x)$.
\end{thm}
\textbf{Proof:} As in the proof of Theorem \ref{thm1}, it can be proved that $\lbrace x_n \rbrace$ is increasing in X and $\lbrace y_n \rbrace$ is decreasing in Y.\\
Now we claim that  \begin{equation}\label{eq15}
d_X(F^{n+1}(x_0,y_0),F^n(x_0,y_0))\leq \Big(\dfrac{b+c}{1-b}\Big)^n\ d_X(x_0,x_1)
\end{equation}
\begin{equation}\label{eq16}
d_Y(G^{n+1}(y_0,x_0),G^n(y_0,x_0))\leq \Big(\frac{a+c}{1-a}\Big)^n\ d_Y(y_0,y_1)
\end{equation}
We prove the claim by mathematical induction, using (\ref{eq13}) and (\ref{eq14}).\\
For $n=1$, consider\\
$d_X(F^2(x_0,y_0),F(x_0,y_0))\\
\hspace*{1.5cm}= d_X(F(F(x_0,y_0),G(y_0,x_0)),F(x_0,y_0))\\
\hspace*{1.5cm}\leq  a\ d_X(F(x_0,y_0),F(x_0,y_0))+b\ d_X(x_0,F^2(x_0,y_0))+ c\ d_X(F(x_0,y_0),x_0)\\
\hspace*{1.5cm}\leq  b\ [ d_X(x_0,F(x_0,y_0))+d_X(F(x_0,y_0),F^2(x_0,y_0))]+ c\ d_X(F(x_0,y_0),x_0)\\
ie,\ d_X(F^2(x_0,y_0),F(x_0,y_0))\leq  \dfrac{b+c}{1-b}\ d_X(x_0,x_1) $\\
Thus the inequality (\ref{eq15}) is true for $n=1$.\\
Now assume that (\ref{eq15}) is true for $n\leq m$, then check for $n=m+1$.\\
Consider,\\
$d_X(F^{m+2}(x_0,y_0),F^{m+1}(x_0,y_0))\\
\hspace*{1.5cm}= d_X(F(F^{m+1}(x_0,y_0),G^{m+1}(y_0,x_0)),F(F^m(x_0,y_0),G^m(y_0,x_0)))\\
\hspace*{1.5cm}\leq  a\ d_X(F^{m+1}(x_0,y_0),F^{m+1}(x_0,y_0))+ b\ d_X(F^m(x_0,y_0),F^{m+2}(x_0,y_0))\\
\hspace*{1.9cm}+c\ d_X(F^{m+1}(x_0,y_0),F^m(x_0,y_0))\\
\hspace*{1.5cm}\leq  b\ [d_X(F^m(x_0,y_0),F^{m+1}(x_0,y_0))+d_X(F^{m+1}(x_0,y_0),F^{m+2}(x_0,y_0))]\\
\hspace*{1.9cm}+c\ d_X(F^{m+1}(x_0,y_0),F^m(x_0,y_0))\\
ie,\ d_X(F^{m+2}(x_0,y_0),F^{m+1}(x_0,y_0))\leq  \dfrac{b+c}{1-b}\ d_X(F^m(x_0,y_0),F^{m+1}(x_0,y_0))\\
\hspace*{6.3cm}\leq \Big(\dfrac{b+c}{1-b}\Big)^{m+1}\ d_X(x_0,x_1)$\\
ie, the inequality (\ref{eq15}) is true for all $n\in \mathbb{N}$\\
Similarly we can prove the inequality (\ref{eq16}).\\
For $m>n$, consider\\
$d_X(F^n(x_0,y_0),F^m(x_0,y_0))\\
\hspace*{1.5cm}\leq  d_X(F^n(x_0,y_0),F^{n+1}(x_0,y_0))+ d_X(F^{n+1}(x_0,y_0),F^{n+2}(x_0,y_0))+...+\\\hspace*{2.5cm}d_X(F^{m-1}(x_0,y_0),F^m(x_0,y_0))\\
\hspace*{1.5cm}\leq  \Bigg[\Big(\dfrac{b+c}{1-b}\Big)^n+\Big(\dfrac{b+c}{1-b}\Big)^{n+1}+...+\Big(\dfrac{b+c}{1-b}\Big)^{m-1}\Bigg]\ d_X(x_0,x_1)\\
\hspace*{1.5cm}\leq  \dfrac{{\delta_1}^n}{1-\delta_1}\ d_X(x_0,x_1);\ \text{where}\ \delta_1= \dfrac{b+c}{1-b}<1\\
\hspace*{1.5cm}\rightarrow 0\  \text{as}\  n\rightarrow \infty $\\
ie, $\lbrace F^n(x_0,y_0)\rbrace$ is a Cauchy sequence in X.\\
Similarly we can prove that $\lbrace G^n(y_0,x_0)\rbrace$ is a Cauchy sequence in Y.\\
Since X and Y are complete, there exist $x\in X$ and $y\in Y$ such that $\lim _{n\rightarrow \infty}F^n(x_0,y_0)=x$ and $\lim _{n\rightarrow \infty}G^n(y_0,x_0)=y$.\\
By continuity of $F$ and $G$, as in the Theorem \ref{thm1} we can show that $F(x,y)=x$ and $G(y,x)=y$.\hspace*{.25cm} $\square$\\

In the following theorem we replace the continuity by other conditions to obtain FG-coupled fixed point.
\begin{thm}\label{thm13}
Let $(X,d_X,\leq_{P_1})$ and $(Y,d_Y,\leq_{P_2})$ be two partially ordered complete metric spaces and $F:X\times Y\rightarrow X$, $G:Y\times X\rightarrow Y$ be two mappings having the mixed monotone property. Assume that X and Y satisfy the following property\begin{enumerate}
\item[(i)] If a non-decreasing sequence $\lbrace x_n \rbrace\rightarrow x$ then $x_n\leq_{P_1} x\ \forall n$
\item[(ii)] If a non-increasing sequence $\lbrace y_n \rbrace\rightarrow y$ then $y\leq_{P_2} y_n\ \forall n$.
\end{enumerate}
Also assume that there exist non-negative a,b,c satisfying\begin{equation}\label{eq30}
\begin{aligned}
d_X(F(x,y),F(u,v))\leq a\ d_X(x,F(u,v))+& b\ d_X(u, F(x,y))+c\ d_X(x,u);\\& \forall x\geq_{P_1}u,\ y\leq_{P_2}v;\ 2b+c<1
\end{aligned}
\end{equation}
\begin{equation}\label{eq31}
\begin{aligned}
d_Y(G(y,x),G(v,u))\leq a\ d_Y(y,G(v,u))+& b\ d_Y(v,G(y,x))+c\ d_Y(y,v);\\& \forall x\leq_{P_1}u,\ y\geq_{P_2}v;\ 2a+c<1
\end{aligned}
\end{equation}
If there exist $x_0\in X, y_0\in Y$ satisfying  $x_0\leq_{P_1} F(x_0, y_0)$ and $y_0\geq_{P_2} G(y_0, x_0)$ then there exist $x\in X, y\in Y$ such that $x=F(x,y)$ and $y=G(y,x)$.
\end{thm}
\textbf{Proof:} Following as in the proof of Theorem \ref{thm5} we get $\lim_{n\rightarrow \infty}F^n(x_0,y_0)=x$ and $\lim_{n\rightarrow \infty}G^n(y_0,x_0)=y$.\\
We have
\begin{align*}
d_X(F(x,y),x)&\leq d_X(F(x,y),F^{n+1}(x_0,y_0))+ d_X(F^{n+1}(x_0,y_0),x)\\
&= d_X(F(x,y),F(F^n(x_0,y_0),G^n(y_0,x_0))+ d_X(F^{n+1}(x_0,y_0),x)\\
&\leq  a\ d_X(x,F(F^n(x_0,y_0),G^n(y_0,x_0)))+ b\ d_X(F^n(x_0,y_0),F(x,y))\\&\ \ \ +c\ d_X(x,F^n(x_0,y_0))+ d_X(F^{n+1}(x_0,y_0),x)\\
&= a\ d_X(x,F^{n+1}(x_0,y_0))+ b\ d_X(F^n(x_0,y_0),F(x,y))\\& \ \ \ +c\ d_X(x,F^n(x_0,y_0))+ d_X(F^{n+1}(x_0,y_0),x)
\end{align*}
ie, $d_X(F(x,y),x)\leq b\ d_X(x,F(x,y))$ as $n\rightarrow \infty$, which implies that $d_X(F(x,y),x)=0$. Therefore $F(x,y)=x$.\\
Also by using (\ref{eq31}) and $\lim_{n\rightarrow \infty}G^n(y_0,x_0)=y$ we can show that $y=G(y,x)$. \hspace*{.25cm} $\square$
\begin{rmrk}\normalfont
If we take $c=0$ in Theorems \ref{thm5} and \ref{thm13}, we get Theorems 2.9 and 2.10 respectively of \cite{8}.
\end{rmrk}
\section{FG-coupled fixed point theorems for generalized quasi-contractions}
The concept of quasi-contraction was defined by Ciric \cite{3} in 1974. A self mapping T on a metric space X is said to be a quasi-contraction iff there exist a number $h$, $0\leq h<1$, such that \begin{equation*}
d(Tx,Ty)\leq h\ \max \lbrace d(x,y),d(x,Tx),d(y,Ty),d(x,Ty),d(y,Tx)\rbrace\end{equation*}  for all $x,y\in X$. In 1979 K.M. Das and K.V. Naik \cite{4} introduced the concept of quasi-contraction for two mappings. Inspired by this we generalize the concept of quasi-contraction to a mapping on product space and prove the following theorems.

\begin{thm}\label{thm6}
Let $(X,d_X,\leq_{P_1})$ and $(Y,d_Y,\leq_{P_2})$ be two partially ordered complete metric spaces. Let $F:X\times Y\rightarrow X$ and $G:Y\times X\rightarrow Y$ be two continuous functions having the mixed monotone property. Assume that there exist $k,l \in \Big[0,\dfrac{1}{2}\Big)$ such that 
\begin{equation}\label{eq17}
d_X(F(x,y),F(u,v))\leq k\ M(x,y,u,v);\ \ \forall x\geq_{P_1}u,\ y\leq_{P_2}v
\end{equation}
\begin{equation}\label{eq18}
d_Y(G(y,x),G(v,u))\leq l\ N(y,x,v,u);\ \ \forall x\leq_{P_1}u,\ y\geq_{P_2}v,
\end{equation} 
where\\
$M(x,y,u,v)= \max \Big\{d_X(x,u),d_X(x,F(x,y)),d_X(x,F(u,v)),d_X(u,F(u,v)),d_X(u,F(x,y))\Big\}$\\
$N(y,x,v,u)= \max \Big\{ d_Y(y,v), d_Y(y,G(y,x)), d_Y(y,G(v,u)), d_Y(v,G(v,u)), d_Y(v,G(y,x)) \Big\}.$\\
If there exist $x_0\in X, y_0\in Y$ satisfying  $x_0\leq_{P_1} F(x_0, y_0)$ and $y_0\geq_{P_2} G(y_0, x_0)$ then there exist $x\in X, y\in Y$ such that $x=F(x,y)$ and $y=G(y,x)$.
\end{thm}
\textbf{Proof:} As in Theorem \ref{thm1}, it can be proved that $\lbrace x_n \rbrace$ is increasing in X and $\lbrace y_n \rbrace$ is decreasing in Y.\\
Now we claim that \begin{equation}\label{eq19}
d_X(F^{n+1}(x_0,y_0),F^n(x_0,y_0))\leq \Big(\dfrac{k}{1-k}\Big)^n\ d_X(x_0,x_1)
\end{equation}
\begin{equation}\label{eq20}
d_Y(G^{n+1}(y_0,x_0),G^n(y_0,x_0))\leq \Big(\frac{l}{1-l}\Big)^n\ d_Y(y_0,y_1)
\end{equation} The proof of the claim is by mathematical induction using (\ref{eq17}) and (\ref{eq18}).\\
For $n=1$, consider\\
$d_X(F^2(x_0,y_0),F(x_0,y_0))\\
\hspace*{.5 cm}= d_X(F(F(x_0,y_0),G(y_0,x_0)),F(x_0,y_0))\\
\hspace*{.5cm}\leq k\ M(F(x_0,y_0),G(y_0,x_0),x_0,y_0)\\
\hspace*{.5cm}= k\ \max \Big\{ d_X(F(x_0,y_0),x_0), d_X(F(x_0,y_0),F^2(x_0,y_0)), d_X(F(x_0,y_0),F(x_0,y_0)),\\ \hspace*{2.5cm} d_X(x_0,F(x_0,y_0)), d_X(x_0,F^2(x_0,y_0))\Big\}\\
\hspace*{.5cm}= k\ \max \Big\{ d_X(x_0,F(x_0,y_0)), d_X(F(x_0,y_0),F^2(x_0,y_0)),  d_X(x_0,F^2(x_0,y_0))\Big\}\\
\hspace*{.5cm}\leq  k\ \max \Big\{ d_X(x_0,F(x_0,y_0)), d_X(F(x_0,y_0),F^2(x_0,y_0)),\\ \hspace*{2.5cm} d_X(x_0,F(x_0,y_0))+ d_X(F(x_0,y_0),F^2(x_0,y_0))\Big\}\\
\hspace*{.5cm}= k\ [d_X(x_0,F(x_0,y_0))+ d_X(F(x_0,y_0),F^2(x_0,y_0))]\\
ie,\ d_X(F^2(x_0,y_0),F(x_0,y_0))\leq  \dfrac{k}{1-k}\ d_X(x_0,F(x_0,y_0))\\
\hspace*{5.2cm}= \dfrac{k}{1-k}\ d_X(x_0,x_1)$\\
So the inequality (\ref{eq19}) is true for $n=1$.\\
Assume that the result is true for $n\leq m$, then check for $n=m+1$.\\
Consider,\\
$d_X(F^{m+2}(x_0,y_0),F^{m+1}(x_0,y_0))\\
= d_X(F(F^{m+1}(x_0,y_0),G^{m+1}(y_0,x_0)),F(F^m(x_0,y_0),G^m(y_0,x_0)))\\
\leq  k\ M(F^{m+1}(x_0,y_0),G^{m+1}(y_0,x_0),F^m(x_0,y_0),G^m(y_0,x_0))\\
=k\ \max \Big\{ d_X(F^{m+1}(x_0,y_0),F^m(x_0,y_0)),d_X(F^{m+1}(x_0,y_0),F^{m+2}(x_0,y_0)),\\ \hspace*{.5cm}d_X(F^{m+1}(x_0,y_0),F^{m+1}(x_0,y_0)),d_X(F^m(x_0,y_0),F^{m+1}(x_0,y_0)),d_X(F^m(x_0,y_0),F^{m+2}(x_0,y_0))\Big\}\\
=k\ \max \Big\{ d_X(F^{m+1}(x_0,y_0),F^m(x_0,y_0)), d_X(F^{m+1}(x_0,y_0),F^{m+2}(x_0,y_0)),\\ \hspace*{1.8cm}d_X(F^m(x_0,y_0),F^{m+2}(x_0,y_0))\Big\}\\
\leq  k\ \max \Big\{ d_X(F^{m+1}(x_0,y_0),F^m(x_0,y_0)), d_X(F^{m+1}(x_0,y_0),F^{m+2}(x_0,y_0)),\\ \hspace*{1.8cm}d_X(F^m(x_0,y_0),F^{m+1}(x_0,y_0))+ d_X(F^{m+1}(x_0,y_0),F^{m+2}(x_0,y_0))\Big\}\\
=k\ [d_X(F^m(x_0,y_0),F^{m+1}(x_0,y_0)) +d_X(F^{m+1}(x_0,y_0),F^{m+2}(x_0,y_0))]\\
ie,\ d_X(F^{m+2}(x_0,y_0),F^{m+1}(x_0,y_0))\leq  \dfrac{k}{1-k}\  d_X(F^{m+1}(x_0,y_0),F^m(x_0,y_0))\\
\hspace*{6.2cm}\leq  \Big(\dfrac{k}{1-k}\Big)^{m+1}d_X(x_0,x_1)$\\
ie, inequality (\ref{eq19}) is true for all $n\in \mathbb{N}$.\\
Similarly we can prove the inequality (\ref{eq20}).\\
Now for $m>n$, consider\\
$d_X(F^m(x_0,y_0),F^n(x_0,y_0))\\
\hspace*{1.5cm}\leq d_X(F^m(x_0,y_0),F^{m-1}(x_0,y_0))+ d_X(F^{m-1}(x_0,y_0),F^{m-2}(x_0,y_0))+...\\ \hspace*{2cm}+d_X(F^{n+1}(x_0,y_0),F^n(x_0,y_0))\\
\hspace*{1.5cm}\leq \Bigg[\Big(\dfrac{k}{1-k}\Big)^{m-1}+\Big(\dfrac{k}{1-k}\Big)^{m-2}+...+\Big(\dfrac{k}{1-k}\Big)^n\Bigg]\ d_X(x_0,x_1)\\
\hspace*{1.5cm}\rightarrow 0\ \text{as}\ n\rightarrow \infty$\\
ie, $\lbrace F^n(x_0,y_0)\rbrace$ is a Cauchy sequence in X.\\
Similarly we can prove that $\lbrace G^n(y_0,x_0)\rbrace$ is a Cauchy sequence in Y.\\
Since X and Y are complete, there exist $x\in X$ and $y\in Y$ such that $\lim _{n\rightarrow \infty}F^n(x_0,y_0)=x$ and $\lim _{n\rightarrow \infty}G^n(y_0,x_0)=y$.\\
As in the Theorem \ref{thm1}, using the continuity of $F$ and $G$ we can show that $x=F(x,y)$ and $y=G(y,x)$. \hspace*{.25cm} $\square$
\begin{eg}\normalfont Let $X= [-1,0],\ Y=[0,1]$ with usual order and usual metric. Define $F: X\times Y\rightarrow X$ by $F(x,y)= \dfrac{x}{3}$ and $G: Y\times X\rightarrow Y$ by $G(y,x)=\dfrac{y}{4}$. Then we can easily check that F satisfies inequality (\ref{eq17}) with $k=\dfrac{1}{3}$ and G satisfies inequality (\ref{eq18}) with $l=\dfrac{1}{4}$ and $(0,0)$ is the FG- coupled fixed point. 
\end{eg}
\begin{thm}\label{thm8}
Let $(X,d_X,\leq_{P_1})$ and $(Y,d_Y, \leq_{P_2})$ be two partially ordered complete metric spaces and $F:X\times Y\rightarrow X$, $G:Y\times X\rightarrow Y$ be two mappings having the mixed monotone property. Assume that X and Y satisfy the following property\begin{enumerate}
\item[(i)]. If a non-decreasing sequence $\lbrace x_n \rbrace\rightarrow x$ then $x_n\leq_{P_1} x\ \forall n$
\item[(ii)]. If a non-increasing sequence $\lbrace y_n \rbrace\rightarrow y$ then $y\leq_{P_2} y_n\ \forall n$.
\end{enumerate}
Also assume that there exist $k,l \in \Big[0,\dfrac{1}{2}\Big)$ such that 
\begin{equation}\label{eq32}
d_X(F(x,y),F(u,v))\leq k\ M(x,y,u,v);\ \ \forall x\geq_{P_1}u,\ y\leq_{P_2}v
\end{equation}
\begin{equation}\label{eq33}
d_Y(G(y,x),G(v,u))\leq l\ N(y,x,v,u);\ \ \forall x\leq_{P_1}u,\ y\geq_{P_2}v
\end{equation} 
where\\
$M(x,y,u,v)= \max \Big\{ d_X(x,u), d_X(x,F(x,y)), d_X(x,F(u,v)), d_X(y,F(u,v)), d_X(u,F(x,y))\Big\}\\
N(y,x,v,u)= \max \Big\{ d_Y(y,v), d_Y(y,G(y,x)), d_Y(y,G(v,u)), d_Y(v,G(v,u)), d_Y(v,G(y,x)) \Big\}$.\\
If there exist $x_0\in X, y_0\in Y$ satisfying  $x_0\leq_{P_1} F(x_0, y_0)$ and $y_0\geq_{P_2} G(y_0, x_0)$ then there exist $x\in X, y\in Y$ such that $x=F(x,y)$ and $y=G(y,x)$.
\end{thm}
\textbf{Proof:}
Following the proof of Theorem \ref{thm11} we get $\lim_{n\rightarrow \infty}F^n(x_0,y_0)=x$ and\\ $\lim_{n\rightarrow \infty}G^n(y_0,x_0)=y$.\\
Now, consider\\ 
$d_X(F(x,y),x)\\
\hspace*{1cm}\leq d_X(F(x,y),F^{n+1}(x_0,y_0))+ d_X(F^{n+1}(x_0,y_0),x)\\
\hspace*{1cm}= d_X(F(x,y),F(F^n(x_0,y_0),G^n(y_0,x_0))+ d_X(F^{n+1}(x_0,y_0),x)\\
\hspace*{1cm}\leq  k\ M(x,y,F^n(x_0,y_0),G^n(y_0,x_0))+ d_X(F^{n+1}(x_0,y_0),x)\\
\hspace*{1cm}= k\ \max \Big\{ d_X(x,F^n(x_0,y_0)),d_X(x,F(x,y)),d_X(x,F^{n+1}(x_0,y_0)),\\\hspace*{1cm} \ \ \ \  d_X(F^n(x_0,y_0),F^{n+1}(x_0,y_0)), d_X(F^n(x_0,y_0),F(x,y))\Big\}+ d_X(F^{n+1}(x_0,y_0),x)$\\
ie, $d_X(F(x,y),x)\leq k\ d_X(x,F(x,y))$ as $n\rightarrow \infty$, which implies that $d_X(F(x,y),x)=0$. Therefore $F(x,y)=x$.\\
Also by using inequality (\ref{eq33}) and  $\lim_{n\rightarrow \infty}G^n(y_0,x_0)=y$, we get  $y=G(y,x)$.\hspace*{.25cm} $\square$
\begin{rmrk}\normalfont
Setting $X=Y$ and $F=G$ in Theorem \ref{thm1} - Theorem \ref{thm8} we get corresponding coupled fixed point theorems in partially ordered complete metric space.
\end{rmrk}
\section*{Acknowledgment}
The first author acknowledges financial support from Kerala State Council for Science, Technology and Environment(KSCSTE), in the form of fellowship.

\hspace*{.1cm}\textsc{Deepa Karichery} (Corresponding author)\\
\hspace*{.1cm}\textsc{Department of Mathematics\\
\hspace*{.1cm}Central University of Kerala, India\\}
\hspace*{.1cm}\textit{E-mail address: deepakarichery@gmail.com}\\

\hspace*{-.5cm}\textsc{Shaini Pulickakunnel\\
\hspace*{.1cm}Department of Mathematics\\
\hspace*{.1cm}Central University of Kerala, India\\}
\hspace*{.1cm}\textit{E-mail address: shainipv@gmail.com}
\end{document}